\documentclass[a4paper, 12pt]{article}
\usepackage{amssymb,latexsym}
\usepackage[cp1251]{inputenc}
\usepackage[english]{babel}
\usepackage{fullpage,amsfonts,amsthm}
\usepackage[intlimits]{amsmath}
\usepackage{indentfirst}
\newtheorem{theorem}{Theorem}
\newtheorem{lemma}[theorem]{Lemma}

\title {
\vspace{15mm}
\textbf{On a diophantine inequality with prime numbers of a special type}}
\author{
\vspace{8mm}
D. I. Tolev %\footnote{Supported by Sofia University, Grant SU 113/2016}
}

\date{}

\begin{document}

\maketitle

\begin{abstract}
We consider the Diophantine inequality 
\[
  \left| p_1^{c} + p_2^{c} + p_3^c- N \right| < (\log N)^{-E} , 
\]
where $1 < c < \frac{15}{14}$,
$N$ is a sufficiently large real number and $E>0$ is an arbitrarily large constant.
We prove that the above inequality has a solution in primes
$p_1$, $p_2$, $p_3$ such that each of the numbers $p_1 + 2, p_2 + 2, p_3 + 2$ has at most
$\left[ \frac{369}{180 - 168 c} \right]$ prime factors, counted with the multiplicity.
\end{abstract}

\section{Introduction and statement of the result}

We consider the diophantine inequality
\begin{equation} \label{100}
 \left| p_1^c + p_2^c + p_3^c - N \right| < \Delta ,
\end{equation}
where $c>1$ is a constant, $N$ is a sufficiently large real number and $\Delta = \Delta(N)$ is a function 
such that $ \Delta(N) \to 0$ as $N \to \infty$. 
Having in mind I.~M.~Vinogradov's famous theorem about Goldbach's ternary problem
(see \cite{Vinogradov}), one may expext that
if $c$ is not much greater than $1$ and $\Delta(N)$ is a suitable function,
then inequality \eqref{100} has a solution in prime numbers $p_1, p_2, p_3$.
A result of this type with $1 < c < \frac{15}{14}$ and with $\Delta = N^{- \kappa}$ for certain
$\kappa = \kappa(c)> 0$
was established in 1992 by the author \cite{Tol_1}.
Several improvements were made since then and the strongest of them is due to Baker and 
Weingartner~\cite{B_W}. In 2014 they established that \eqref{100} is solvable in primes,
provided that $N$ is large enough, 
$1 < c < \frac{10}{9}$ and $\Delta = N^{- \kappa}$ for certain $\kappa = \kappa(c)> 0$.

\bigskip

Suppose that $r$ is a natura number and let $\mathcal P_r$ be the set of positive integers 
having at most $r$ prime factors, counted with the multiplicity. 
(We say that the numbers from $\mathcal P_r$ are
almost primes of order $r$.)
In 1973 Chen~\cite{Chen}, improving results of other mathematicians, established 
that there exist infinitely many primes $p$ such that
$p+2 \in \mathcal P_2$. 
Bearing in mind Chen's result, one may try to study the arithmetical properties
of the set of primes $p$ such that $p+2 \in \mathcal P_r$ for a fixed $r \ge 2$
and, in particular, to establish the solvability of diophantine equations or inequalities 
in such primes.
For example, Matom\"aki and Shao~\cite{Matomaki_Shao}, improving author's results from \cite{Tol_2} and \cite{Tol_3}
as well as a result of  Matom\"aki~\cite{Matomaki},
proved that every sufficiently large odd integer $N$ can be represented as a sum ot three
primes $p_1, p_2, p_3$ such that $p_i + 2 \in \mathcal P_2$, $i=1,2,3$.
Other results of this type were found by the author~\cite{Tol_4}, and by
Dimitrov and Todorova~\cite{Dim_Tod}.

\bigskip

One may expect that 
if the constant $c>1$ is close to one then 
inequality \eqref{100},
with a suitable $\Delta$ satisfying $\Delta \to 0 $ as $N \to \infty$,
is solvable in primes $p_i$ such that $p_i + 2 $ are almost primes of certain fixed order. 
An attempt to establish a result of this type was made by Dimitrov~\cite{Dim_Thesis}, 
but he consideres the inequality \eqref{100} only when $c< \frac{4}{21}$, whilst
the case $c>1$ is more interesting.
Dimitrov recently announsed that he is able to 
prove the solvability of \eqref{100} in the case $1 < c < \frac{121}{120}$ but such a result 
has not been published.

\bigskip

In the present paper we 
assume that
$c$ is a constant such that
\begin{equation} \label{101}
  1 < c < \frac{15}{14} .
\end{equation} 
We consider the inequality
\eqref{100} with 
\begin{equation} \label{102.5}
  \Delta = (\log N)^{-E} ,
\end{equation} 
where $E>0$ is an arbitrarily large constant and 
we prove the following
\begin{theorem} \label{theorem10}
Let $c$ be a constant satisfying \eqref{101}
and let $N$ be a sufficiently large real number.
Then inequality \eqref{100}, with $\Delta$ specified by 
\eqref{102.5}, 
has a solution in primes $p_1, p_2, p_3$ such that each of the numbers $p_i + 2$, $i = 1,2,3$
has at most $\left[ \frac{369}{180 - 168 c} \right]$ prime factors, counted with the multiplicity.
\end{theorem}

It follows from Theorem~\ref{theorem10}
that if $c>1$ is close to $1$, then inequality \eqref{100}, with $\Delta$ given by \eqref{102.5}, 
has a solution in primes $p_i$ such that $p_i + 2 \in \mathcal P_{30}$.

\bigskip

We can establish a similar result for the inequality \eqref{100} with $\Delta = N^{- \kappa}$ for certain
$\kappa > 0$, but then we would have $p_i + 2 \in \mathcal P_m$, $i=1,2,3$, where
$m$ depends on $c$ and $\kappa$.

\bigskip

Notations in the paper shall be as follows.
By $\varepsilon$ and $A$ we denote an
arbitrarily small positive number and respectively, an arbitrarily large constant
which may not be the same in different formulae.
The letter $p$ always denotes a prime number. 
By $\tau(n)$, $\mu(n)$, $\varphi(n)$ and $\Lambda(n)$ we denote 
the number of divisors of $n$,
M\"obius' function, Euler's function
and Von Mangoldt's function respectively.
We shall use $(m, n)$ and $[m, n]$ for the greatest common divisor and the least common multiple of
the integers $m, n$. (We denote in this way
also open and closed intervals from the real line, 
but the meaning will be clear from the context).
Let $[t]$ be the integer part of the real number $t$ and $e(t) = e^{2 \pi i t}$.
With $\chi$ we denote a Dirichlet's character. As usual, $\sum_{\chi \pmod{q}}$ means that the summation
is taken over all Dirichlet's characters modulo $q$. Respectively, $\sum_{\chi \pmod{q}^*}$
means that the summation is taken over the primitive Dirichlet's characters modulo $q$.

\bigskip

Suppose that $\chi$ is a Dirichlet's character and $L(s, \chi)$ is the corresponding $L$-function.
If $T \ge 2$ and $0 \le \sigma \le 1$
we denote by $N(T, \sigma, \chi)$ the number of zeros $\rho = \beta + i \gamma $ of $L(s, \chi)$ 
such that $|\gamma| \le T$ and $\sigma \le \beta \le 1$.
We also write $N(T, \chi) = N (T, 0, \chi)$.

\bigskip

By
\begin{equation} \label{50}
  \psi(y) = \sum_{n \le y} \Lambda(n) , \qquad
  \psi (y, k, l) = \sum_{\substack{n \le y \\ n \equiv l \pmod{k}}} \Lambda(n)
\end{equation}
we denote Chebyshev's functions and we
define $\Delta(y, k, l)$ by
\begin{equation} \label{55}
  \Delta(y, k, l) = \psi (y, k, l) - \frac{y}{\varphi(k)} .
\end{equation}
For a given Dirichlet's character $\chi$ we write
\begin{equation} \label{55.1}
  \psi(y, \chi) = \sum_{n \le y} \Lambda(n) \chi(n) .
\end{equation}

\bigskip

Finally, by $\square$ we mark an end of a proof or its absence.

\section{Beginning of the proof}

Let
$\eta$, $\delta$, $\xi$, $\mu$ be positive real numbers depending on $c$.
We shall specify them later 
but for now only assume that
they satisfy the conditions
\begin{equation} \label{102}
  \xi + 3 \delta < \frac{12}{25} , \qquad 2 < \frac{\delta}{\eta} < 3 , \qquad \mu < 1 .
\end{equation} 
We define
\begin{equation} \label{103}
    X = N^{\frac{1}{c}} , \qquad z = X^{\eta} \qquad D = X^{\delta} ,
  \qquad \tau = X^{\xi - c} ,
\end{equation} 

\begin{equation} \label{103.5}
      r = \left[ (\log X)^2 \right] ,
     \qquad \Xi = (\log X)^{E+3} 
\end{equation} 
and
\begin{equation} \label{104}
  P(z) = \prod_{2 < p < z} p ,
\end{equation} 
where the product is taken over prime numbers.

\bigskip

Consider the sum
\begin{equation} \label{150}
  \Gamma = \sum_{\substack{\mu X < p_1, p_2, p_3 \le X \\ \eqref{100}, \; \eqref{120} }}
   (\log p_1) (\log p_2) (\log p_3) ,
\end{equation}
where the summation is taken over the primes $p_1, p_2, p_3$ from the interval $(\mu X, X]$ which satisfy 
\eqref{100} (with $\Delta$ given by \eqref{102.5}),
as well as the conditions 
\begin{equation} \label{120}
  (p_1 + 2 , P(z)) = (p_2 + 2 , P(z)) = (p_2 + 2 , P(z)) = 1 .
\end{equation}

\bigskip

If we prove the inequality
\begin{equation} \label{119}
   \Gamma > 0 ,
\end{equation}
then the equation \eqref{100} would have a solution 
in primes $p_1, p_2, p_3$ satisfying \eqref{120}. If the number $p_i + 2$ has $l$ prime factors
counted with multiplicity, then from  \eqref{103}, \eqref{104}, \eqref{120} and from the condition 
$\mu X < p_i \le X$ we easily find that 
\begin{equation} \label{120.5}
  l \le \frac{1}{\eta} .
\end{equation}
This means that $p_i + 2$ would be an almost-prime of order $\left[ \eta^{-1} \right]$.
Therefore, to prove the theorem we have to establish \eqref{119} for a suitable choice of $\eta$.

\bigskip

Firstly, we use the following
\begin{lemma} \label{lemma1}
Let $a, \delta$ be real numbers, $0 < \delta < \frac{a}{4}$, and let $r$ be a positive integer.
There exists a function $\theta(y)$ which is $r$ times continuously differentiable and
such that
\begin{align}
 \theta(y) 
   & =
    1 \qquad \text{for} \quad |y| \le a - \delta , 
     \notag  \\
  0 < \theta(y) & < 1  \qquad \text{for} \quad a - \delta < |y| \le a + \delta , 
   \notag \\
 \theta(y) 
   & =
    0 \qquad \text{for} \quad |y| \ge a + \delta  ,  
    \notag 
\end{align}
and its Fourier transform
\begin{equation} \label{157}
   \Theta(x) = \int_{- \infty}^{\infty} \theta(y) e(-xy) \, d y
\end{equation}
satisfies the inequality
\begin{equation} \label{158}
   |\Theta(x)| \le \min \left( 2a , \frac{1}{|x|}   \left( \frac{r}{|x| \delta} \right)^r \; \right) .
\end{equation}
\end{lemma}

\bigskip

{\bf Proof.} This is an old result which goes back to Segal \cite{Segal}.

\hfill
$\square$

\bigskip

We apply this lemma with $r$ given by \eqref{103.5} 
and with $a = \frac{7 \Delta}{8}$, $\delta = \frac{\Delta}{8}$, where $\Delta$
is specified by \eqref{102.5}. Hence we have
\begin{equation} \label{155}
   \theta(y) = 0 \quad \text{for} \quad |y| \ge \Delta ,
   \qquad
   0 < \theta(y) \le 1 \quad \text{for} \quad   |y| <  \Delta 
\end{equation}
and from \eqref{150} and \eqref{155} we find that
\begin{equation} \label{158.1}
   \Gamma \ge \Gamma' := 
   \sum_{\substack{\mu X < p_1, p_2, p_3 \le X \\  \eqref{120} }}
   (\log p_1) (\log p_2) (\log p_3) \, \theta \left( p_1^c + p_2^c + p_3^c - N \right)
\end{equation}

\bigskip

For $i=1, 2, 3$ we consider the quantities
\begin{equation} \label{158.2}
  \Lambda_i = \sum_{d \mid  (p_i + 2 , P(z)) } \mu(d) 
 =   \begin{cases}
   1  \quad & \text{if} \; (p_i + 2 , P(z)) = 1 , \\
   0 & \text{otherwise} .
  \end{cases}
\end{equation}
Bearing in mind  \eqref{120} and \eqref{158.1}
we find that
\begin{equation} \label{160}
  \Gamma' = \sum_{\mu X < p_1, p_2, p_3 \le X }
   (\log p_1) (\log p_2) (\log p_3) \; \Lambda_1 \, \Lambda_2 \, \Lambda_3 \;
    \theta \left( p_1^c + p_2^c + p_3^c - N \right) .
\end{equation}

\bigskip

Now we apply the following fundamental result from sieve theory
\begin{lemma} \label{lemma2}
Suppose that $D > 4$ is a real number.
There exist arithmetical functions $\lambda^{\pm}(d)$
(called Rosser's functions of level $D$) with the following properties.

\bigskip

1) 
For any positive integer $d$ we have
\begin{equation} \label{165}
  |\lambda^{\pm}(d)| \le 1 , \qquad
  \lambda^{\pm}(d) = 0 \quad \text{if} \quad  d > D \quad \text{or} \quad \mu(d) = 0 .
\end{equation}

2) 
If $n$ is a positive integer then 
\[
  \sum_{d \mid n} \lambda^{-} (d) \le \sum_{d \mid n} \mu(d) \le \sum_{d \mid n} \lambda^{+} (d)  .
\]

3) 
If $z$ is a real number such that $z^2 \le D \le z^3$ and if 
\begin{equation} \label{184}
 P(z) = \prod_{2 < p < z} p \, , \quad 
  \mathfrak P = \prod_{2 < p < z} \left( 1 - \frac{1}{p-1} \right) \, , \quad 
   \mathfrak  N^{\pm} = \sum_{ d \mid P(z)} \frac{\lambda^{\pm}(d)}{\varphi(d)} ,
   \quad
  s_0 = \frac{\log D}{\log z} ,
\end{equation}
then we have
\begin{align} 
    \mathfrak P 
  \le
  \mathfrak N^{+}
  \le 
  \mathfrak P 
  \left( F(s_0) + O \left( (\log D)^{- \frac{1}{3}} \right)  \right) ,
   \label{181} \\
   & \notag \\
    \mathfrak P 
  \ge
  \mathfrak N^{-}
  \ge 
  \mathfrak P 
  \left( f(s_0) + O \left( (\log D)^{- \frac{1}{3}} \right)  \right) ,
   \label{182} 
\end{align}
where
$F(s)$ and $f(s)$ (the functions of the linear sieve) satisfy
\begin{equation} \label{185}
 f(s) = 2 e^{\gamma} s^{-1} \log (s-1) , \qquad
 F(s) = 2 e^{\gamma} s^{-1} \qquad \text{for} \qquad 2 \le s \le 3 .
\end{equation}
Here $\gamma$ stands for the Euler constant.
\end{lemma}

\bigskip

{\bf Proof.}
This is a special case of a more general result --- see Greaves \cite[Ch.~4]{Greaves}.

\hfill
$\square$

\bigskip

Suppose that $\lambda^{\pm}(d)$ are the Rosser functions of level $D$, where $D$ is defined by \eqref{103}.
According to Lemma~\ref{lemma2}, if we denote
\begin{equation} \label{170}
   \Lambda_i^{\pm} = \sum_{d \mid  (p_i + 2 , P(z)) } \lambda^{\pm}(d) , \qquad i = 1, 2, 3
\end{equation}
and if $\Lambda_i$ is given by \eqref{158.2},
then we have
\begin{equation} \label{180}
  \Lambda_i^{-} \le \Lambda_i \le \Lambda_i^{+}  , \qquad \Lambda_i \in \{0, 1\} .
\end{equation}

\bigskip

Next we use the following
\begin{lemma} \label{lemma3}
Suppose that $\Lambda_i, \Lambda_i^{\pm}$, $i = 1,2,3$ are real numbers satisfying \eqref{180}.
Then we have
\begin{equation} \label{190}
  \Lambda_1 \, \Lambda_2 \, \Lambda_3 \ge 
  \Lambda_1^{-} \, \Lambda_2^{+} \, \Lambda_3^{+} 
  + \Lambda_1^{+} \, \Lambda_2^{-} \, \Lambda_3^{+} 
  + \Lambda_1^{+} \, \Lambda_2^{+} \, \Lambda_3^{-}
  - 2 \Lambda_1^{+} \, \Lambda_2^{+} \, \Lambda_3^{+} .
\end{equation}
\end{lemma}

\bigskip

{\bf Proof.}
The proof is elementary and similar to the proof of
\cite[Lemma~13]{Br_F}).

\hfill
$\square$

\bigskip

We apply \eqref{160}, \eqref{180} and \eqref{190}
and then we substitute the quantity from the right side of \eqref{190}
for $\Lambda_1 \, \Lambda_2 \, \Lambda_3 $ in \eqref{160}. We find that
\[
  \Gamma' \ge \Gamma_1 + \Gamma_2 + \Gamma_3 - 2 \Gamma_4 ,
\]
where $\Gamma_1, \dots, \Gamma_4$ are the contributions comming from 
the consequtive terms of the right side of \eqref{190}.
It is clear that
\begin{align}
 \Gamma_1 = \Gamma_2 = \Gamma_3 
  & =
  \sum_{\mu X < p_1, p_2, p_3 \le X }
   (\log p_1) (\log p_2) (\log p_3) \; \Lambda_1^{-} \, \Lambda_2^{+} \, \Lambda_3^{+} \;
    \theta \left( p_1^c + p_2^c + p_3^c - N \right) , 
   \label{200} \\
   & \notag \\
   \Gamma_4
   & =
 \sum_{\mu X < p_1, p_2, p_3 \le X }
   (\log p_1) (\log p_2) (\log p_3) \; \Lambda_1^{+} \, \Lambda_2^{+} \, \Lambda_3^{+} \;
   \theta \left( p_1^c + p_2^c + p_3^c - N \right) .
   \label{210}
\end{align}
Hence, we get
\begin{equation} \label{212}
  \Gamma' \ge 3 \Gamma_1 - 2 \Gamma_4 .
\end{equation}

\bigskip

Consider $\Gamma_1$. (The study of $\Gamma_4 $ is likewise).
We apply Fourier's inversion formula
\[
  \theta(t) = 
  \int_{- \infty}^{\infty} \Theta(x) \, e(xt) \, d x
\]
as well as \eqref{170} to find that
\begin{align}
  \Gamma_1 
  & = 
   \sum_{\mu X < p_1, p_2, p_3 \le X }
   (\log p_1) (\log p_2) (\log p_3) \; \Lambda_1^{-} \, \Lambda_2^{+} \, \Lambda_3^{+} \;
   \int_{- \infty}^{\infty} \Theta(x) e \left(x \left( p_1^c + p_2^c + p_3^c - N \right) \right) \, d x
   \notag \\
   & \notag \\
   & = 
   \int_{- \infty}^{\infty} \Theta(x) \, L^{-}(x) \, L^{+}(x)^2 \, e (- N x) \, d x ,
   \label{230}
\end{align}
where
\[
  L^{\pm}(x) = \sum_{\mu X < p \le X} (\log p) \, 
  e \left( x p^c \right) \, 
  \sum_{d \mid (p+2, P(z))} \lambda^{\pm} (d) .
\]
Changing the order of summation, we get
\begin{equation} \label{240}
  L^{\pm}(x) = \sum_{d \mid P(z)} \lambda^{\pm}(d) \sum_{\substack{\mu X < p \le X \\ d \mid p+2 }} (\log p) \, 
  e \left( x p^c \right) .
\end{equation}

\bigskip

We divide the integral from \eqref{230} into three parts as follows:
\begin{equation} \label{250}
  \Gamma_1  = \Gamma_1^{(1)} + \Gamma_1^{(2)} + \Gamma_1^{(3)} ,
\end{equation}
where
\begin{align}
  \Gamma_1^{(1)} 
  & = \int_{|x| < \tau} \Theta(x) \, L^{-}(x) \, L^{+}(x)^2 \, e (- N x) \, d x ,
  \label{260} \\
  & \notag \\
  \Gamma_1^{(2)} 
  & = \int_{\tau < |x| <  \Xi } \Theta(x) \, L^{-}(x) \, L^{+}(x)^2 \, e (- N x) \, d x ,
  \label{270} \\
  & \notag \\
  \Gamma_1^{(3)} 
  & = \int_{|x| >  \Xi } \Theta(x) \, L^{-}(x) \, L^{+}(x)^2 \, e (- N x) \, d x .
  \label{280} 
\end{align}

\bigskip

Similarly, for the quantity $\Gamma_4$ defined by \eqref{210} we find
\begin{equation} \label{250.1}
  \Gamma_4  = \Gamma_4^{(1)} + \Gamma_4^{(2)} + \Gamma_4^{(3)} ,
\end{equation}
where
\begin{align}
  \Gamma_4^{(1)} 
  & = \int_{|x| < \tau} \Theta(x) \, L^{+}(x)^3 \, e (- N x) \, d x ,
  \label{260.1} \\
  & \notag \\
  \Gamma_4^{(2)} 
  & = \int_{\tau < |x| <  \Xi } \Theta(x)  \, L^{+}(x)^3 \, e (- N x) \, d x ,
  \label{270.1} \\
  & \notag \\
  \Gamma_4^{(3)} 
  & = \int_{|x| >  \Xi } \Theta(x)  \, L^{+}(x)^3 \, e (- N x) \, d x .
  \label{280.1} 
\end{align}

\bigskip

It is easy to estimate $\Gamma_1^{(3)}$ and $\Gamma_4^{(3)}$.
It is clear from \eqref{240} that $L^{\pm}(x) \ll X^{1 + \varepsilon}$. 
We also use \eqref{102.5}, \eqref{103.5} and \eqref{158} to find that
\begin{equation} \label{287.5}
  \Gamma_1^{(3)} , \Gamma_4^{(3)} \ll 
  X^{3 + \varepsilon} \left( \frac{r}{\Delta} \right)^{r} \int_{ \Xi }^{\infty} \frac{d x}{x^{r+1}}
  \ll X^{3+ \varepsilon} \left( \frac{r}{\Delta  \, \Xi} \right)^{r} \ll 1 .
\end{equation}

\bigskip

From \eqref{158.1}, \eqref{212}, \eqref{250} and \eqref{250.1} it follows that
\begin{equation} \label{288.5}
  \Gamma \ge 
  \left|  3 \Gamma_1^{(1)} - 2 \Gamma_4^{(1)} \right|
   - c_0 \left( 
   \big|  \Gamma_1^{(2)} \big|
  + \big|  \Gamma_4^{(2)} \big|  + 1 \right)  ,
\end{equation}
where $c_0 > 0$ is an absolute constant.

\section{The integrals $\Gamma_1^{(1)}$ and $\Gamma_4^{(1)}$}

In this section we find asymptotic formulae for $L^{\pm}(x)$, provided that $|x|< \tau$.
The arithmetic structure of the Rosser weights $\lambda^{\pm}(d)$ is not important here, so
we consider a sum of the form
\begin{equation}  \label{289} 
  L(x) 
   = 
  \sum_{d \le D} \lambda(d) \; \sum_{\substack{ \mu X < p \le X \\ d \mid p+2}}
  (\log p) \, e \left( x p^c \right) ,
\end{equation}
and we assume that $\lambda(d)$ are real numbers satisfying
\begin{equation} \label{288}
 |\lambda(d)| \le 1 , \qquad
 \lambda(d) = 0 \quad \text{if} \quad 2 \mid d \quad \text{or} \quad \mu(d) = 0 .
\end{equation}
We also define
\begin{equation} \label{290}
  I(x)
   =
    \int_{\mu X}^X e \left( x t^c \right) \, d t ,
\end{equation}

\bigskip

To study $I(x)$, we need the following 
\begin{lemma} \label{lemma3.5}
Consider the integral 
\[
  I = \int_a^b G(x) \, e (F(x)) \, d x
\]
where
$G(x), F(x)$ are real functions with continuous second derivatives.
Assume that the function $G(x)/F'(x)$ is monotonous
and suppose that $|G(x)| \le H$ for all $x \in [a, b]$.

\bigskip

If $|F'(x)| \ge  h > 0$ for all $x \in [a, b]$ 
then $ I \ll H h^{-1}  $.

\bigskip

If $|F''(x)| \ge h > 0$ for all $ x \in [a, b]$ then
$  I \ll H h^{- \frac{1}{2}} $.
\end{lemma}

{\bf Proof.} See \cite[p. 71]{Titchmarsh}.

\hfill
$\square$

\bigskip

We also need Bombieri--Vinogradov's theorem:
\begin{lemma} \label{lemma4}
Suppose that $x > 2$, $Q > 2$ and consider the sum
\[
  \Sigma = \sum_{d \le Q} \max_{y \le x} \max_{(l, d)=1} \left| \Delta(y, d, l) \right| ,
\]
where $\Delta(y, d, l)$ is defined by \eqref{55}. 
For any constant $A>0$ there exists $B=B(A)>0$ such that if $Q \le \sqrt{x} \, (\log x)^{-B}$,
then $ \Sigma \ll x (\log x)^{-A} $.
\end{lemma}

\bigskip

{\bf Proof.}
See \cite[Ch. 28]{Davenport}.

\hfill
$\square$

\bigskip

In the next lemma we give explicit formulae for Chebyshev's function $\psi(y)$
and for the function $\psi(y, \chi)$ defined by \eqref{55.1}.
\begin{lemma} \label{lemma4.5}
Suppose that $2 \le T \le y$.

\bigskip

We have
\begin{equation} \label{449}
 \psi(y) = y - \sum_{|\gamma| < T} \frac{y^{\rho}}{\rho} + O \left( \frac{y \log y}{T} \right) , 
\end{equation}
where the summation runs over the non-trivial zeros $\rho = \beta + i \gamma $ 
of the Riemann zeta function such that $|\gamma| < T $.

\bigskip

If $r>1$ and $\chi$ is a primitive character $\pmod{r}$, then
we have
\begin{equation} \label{450}
  \psi(y, \chi) = - \sum_{|\gamma| < T} \frac{y^{\rho}}{\rho} 
  + \sum_{|\gamma| < 1} \frac{1}{\rho} + O \left( \frac{y \log (ry)}{T} \right) ,
\end{equation}
where the summation is taken over the non-trivial zeros $\rho = \beta + i \gamma$ of Dirichlet's 
$L$-function $L(s, \chi)$ such that $|\gamma| < T$ and $|\gamma| < 1$, respectively.
\end{lemma}

\bigskip

{\bf Proof.} See \cite[Ch.~17]{Davenport} and
\cite[Ch.~19]{Davenport}

\hfill
$\square$

\bigskip

The next lemma provides an information about the density of the zeroes of
Dirichlet's $L$-functions.
\begin{lemma} \label{lemma4.6}
If $\chi$ is a primitive Dirichlet's chracter modulo $d$ and if $T \ge 2$, then we have
\begin{equation} \label{450.1}
  N(T, \chi) \ll T \log (dT) .
\end{equation}
If $Q \ge 1$ and $T \ge 2$, then for the sum
\begin{equation} \label{450.2}
  \Sigma (T, \sigma, Q) = \sum_{d \le Q} \; \sum_{\chi \!\!\!\! \pmod{d}^*} N(T, \sigma, \chi) 
\end{equation}
we have
\begin{equation} \label{450.3}
  \Sigma (T, \sigma, Q) \ll
   \begin{cases}
   (Q^2 T)^{\frac{3(1 - \sigma)}{2 - \sigma}} \left( \log (QT) \right)^9  
      \quad \text{if} \quad \frac{1}{2} \le \sigma \le \frac{4}{5} ,   \\
       \\
   (Q^2 T)^{\frac{2(1 - \sigma)}{\sigma}} \left( \log (QT) \right)^{14}  
      \quad \text{if} \quad \frac{4}{5} \le \sigma \le 1 .
   \end{cases}
\end{equation}
\end{lemma}

{\bf Proof.}
The proof of \eqref{450.1} can be found in 
\cite[Ch.~16]{Davenport} and for the proof of \eqref{450.3} see
\cite[Theorem~12.2]{Montgomery}. 

\hfill
$\square$

\bigskip

In the next lemma we present an analog of the estimate from
\cite[Lemma~5]{Tol_1}. 
\begin{lemma} \label{lemma5}
Consider the sum
\begin{equation} \label{660}
 \mathfrak L (T, Q, X) = \sum_{d \le Q} \; \sum_{\chi \!\!\!\! \pmod{d}^* } \; \sum_{|\gamma| \le T} X^{\beta} ,
\end{equation}
where $X \ge 2$, $Q \ge 1$, $T \ge 2$ and where the summation in the inner sum is taken over 
the non-trivial zeros $\rho = \beta + i \gamma$ of Dirichlet's $L$-function $L(s, \chi)$
such that $|\gamma| \le T$.

\bigskip

Suppose that
\begin{equation} \label{679.5}
  Q^2 T \le X^{\frac{12}{25}} .
\end{equation}
Then we have
\begin{equation} \label{680}
  \mathfrak L (T, Q, X) \ll (\log X)^{20} \left( X + X^{\frac{4}{5}} \, Q \, T^{\frac{1}{2}} \right) .
\end{equation}
\end{lemma}

{\bf Proof.}
The proof of \eqref{680} is standard but for reader's convenience we present the arguments.
Suppose that $\chi$ is a 
primitive charcacter
modulo $d$ and 
$\rho = \beta + i \gamma$ is a non-trivial zero of $L(s, \chi)$.
We start with the identity
\begin{equation} \label{680.1}
 X^{\beta} = 1 + (\log X) \int_0^1 X^{\sigma} \, \Phi(\sigma, \beta) \, d \sigma ,
\end{equation}
where $\Phi(\sigma, \beta) = 0 $ for $\sigma > \beta$ and 
$\Phi(\sigma, \beta) = 1 $ for $ \sigma \le \beta$.
It is clear that
\[
  \sum_{|\gamma| \le T} \Phi(\sigma, \beta) = N(T, \sigma, \chi) ,
\]
hence applying \eqref{680.1} and the estimate \eqref{450.1} from Lemma~\ref{lemma4.6} we find that
\begin{align}
  \sum_{|\gamma| \le T} X^{\beta} 
   & = 
  N(T, \chi) + (\log X) \int_0^1 X^{\sigma} \, N(T, \sigma, \chi) \, d \sigma 
   \notag \\
   & \notag \\
   & \ll 
    X^{\frac{1}{2}} T (\log X)^2 +
    (\log X) \int_{\frac{1}{2}}^1 X^{\sigma} \, N (T, \sigma, \chi) \, d \sigma .
  \notag 
\end{align}
From the above formula, \eqref{450.2} and \eqref{660} we find that
\[
  \mathfrak L (T, Q, X) \ll
  X^{\frac{1}{2}} Q^2 T (\log X)^2 +
  (\log X) \int_{\frac{1}{2}}^1 X^{\sigma} \, \Sigma (T, \sigma, Q) \, d \sigma .
\]
We apply the estimate \eqref{450.3} from Lemma~\ref{lemma4.6} and find that
\begin{equation} \label{680.2}
  \mathfrak L (T, Q, X) \ll
  (\log X)^{15} \left( X^{\frac{1}{2}} Q^2 T  + \mathcal I_1 + \mathcal I_2 \right) ,
\end{equation}
where
\[
  \mathcal I_1 = \int_{\frac{1}{2}}^{\frac{4}{5}} X^{\sigma} \, (Q^2 T)^{\frac{3(1-\sigma)}{2 - \sigma}}
   \, d \sigma ,
  \qquad
  \mathcal I_2 = \int_{\frac{4}{5}}^1 X^{\sigma} \, (Q^2 T)^{\frac{2(1-\sigma)}{\sigma}} \, d \sigma .
\]
To estimate the integral $\mathcal I_1$, we write it in the form
\[
  \mathcal I_1 = \int_{\frac{1}{2}}^{\frac{4}{5}}  e^{ h_1(\sigma) } \, d \sigma , \qquad
  h_1(\sigma) = \sigma (\log X) + \frac{3(1-\sigma)}{2 - \sigma} \log(Q^2 T) .
\]
Using the condition \eqref{679.5},
it is easy to verify that $h'_1(\sigma) \ge 0$ for $\frac{1}{2} \le \sigma \le \frac{4}{5}$.
This means that $\max_{\frac{1}{2} \le \sigma \le \frac{4}{5}} h_1(\sigma) = h_1 \left( \frac{4}{5} \right)$
and therefore
\begin{equation} \label{680.3}
 \mathcal I_1  \ll e^{ h_1 \left( \frac{4}{5} \right) } = 
 X^{\frac{4}{5}} Q T^{\frac{1}{2}} .
\end{equation}

\bigskip

Consider $\mathcal I_2 $. We have
\[
  \mathcal I_2 = \int_{\frac{4}{5}}^1  e^{ h_2(\sigma) } \, d \sigma , \qquad
  h_2(\sigma) =  \sigma (\log X)  + \frac{2(1 - \sigma)}{\sigma} \log(Q^2 T) .
\]
We have
$h''_2(\sigma) > 0 $ for $\frac{4}{5} \le \sigma \le 1$, hence
$ \max_{\frac{4}{5} \le \sigma \le 1} h_2(\sigma) = 
   \max \left( h_2 \left( \frac{4}{5} \right) , h_2 \left( 1 \right)\right) $ 
and therefore
\begin{equation} \label{680.4}
  \mathcal I_2  \ll e^{h_2 \left( \frac{4}{5} \right)} + e^{h_2 \left( 1 \right)} =
  X^{\frac{4}{5}} Q T^{\frac{1}{2}} + X .
\end{equation}

\bigskip

The estimate \eqref{680} is a consequence of \eqref{679.5} and \eqref{680.2} -- \eqref{680.4}.

\hfill
$\square$

\bigskip

We shall prove the following 
\begin{lemma} \label{lemma10}
Suppose that $D$ and $\tau$ are defined by \eqref{103} and that $\xi$ and $\delta$ satisfy
\eqref{102}. 
If $L(x)$ and $I(x)$ are defined by \eqref{289} and \eqref{290}
and if $|x| < \tau $ then we have
\begin{equation}  \label{300}
  L(x) = \sum_{ d \le D} \frac{\lambda(d)}{\varphi(d)} \; I(x) 
  + O \left( \, X (\log X)^{-A} \, \right) ,
\end{equation}
where $A>0$ is an arbitrarily large constant.
\end{lemma}

{\bf Proof.}
One may easily see that
\begin{equation}  \label{330}
  L(x) = L_1 (x) + O \left( X^{\frac{1}{2} + \varepsilon }\right) ,
\end{equation}
where
\begin{equation}  \label{340}
  L_1 (x) =
  \sum_{d \le D} \lambda(d) \, S(x, d) , \qquad
  S(x, d) =   
  \sum_{\substack{\mu X < n \le X \\ d \mid n+2 }} \Lambda(n) \, 
  e \left( x n^c \right) .
\end{equation}
To find an asymptotic formula for $L_1(x)$, we shall proceed in different ways 
according to the size of $|x|$.

\bigskip

Firstly, consider the case
\begin{equation}  \label{350}
  |x| \le X^{-c} (\log X)^B ,
\end{equation}
where $B>0$ is a constant which we shall specify later.
For the sum $S(x, d)$, defined by \eqref{340},
we apply Abel's formula and \eqref{50} to find that
\begin{align}
 S(x, d) 
   & = 
   e \left( x X^c \right)  \sum_{\substack{\mu X < n \le X \\ d \mid n+2 }} \Lambda(n) 
   - \int_{\mu X}^X
    \sum_{\substack{\mu X < n \le t \\ d \mid n+2 }} \Lambda(n) \, 
    \frac{d}{d t} e \left( x t^c \right) \, d t 
    \notag \\
    & \notag \\
   & = 
  e \left( x X^c \right) \big(  \psi(X, d, -2) - \psi(\mu X, d, -2) \big)
   \notag \\
    & \notag \\
   & \qquad \qquad \qquad \qquad \qquad 
  - \int_{\mu X}^X \big( \psi(t, d, -2) - \psi(\mu X, d, -2) \big) \, \frac{d}{d t} e \left( x t^c \right) .
  \notag
\end{align}
Now we use \eqref{55} to get
\begin{align}
  S(x, d) 
   & = 
   e \left( x X^c \right) 
   \left( \frac{X - \mu X}{\varphi(d)} + \Delta(X, d, -2) - \Delta(\mu X, d, -2) \right)
   \notag \\
   & \notag \\
   & \qquad \qquad \qquad
   - \int_{\mu X}^X \left( \frac{t - \mu X}{\varphi(d)} + \Delta(t, d, -2) - \Delta(\mu X, d, -2) \right) 
   \frac{d}{d t} e \left( x t^c\right) \, d t .
   \notag
\end{align}
Therefore, using \eqref{290} and the assumption \eqref{350} we find that
\begin{align}
  S(x, d) 
   & = 
   \frac{1}{\varphi(d)} 
    \left( e \left( x X^c\right) (X - \mu X) - \int_{\mu X}^X
    (t - \mu X) \, \frac{d}{d t} e \left( x t^c \right) 
    \right) 
     \notag \\
     & \notag \\
     & \qquad \qquad \qquad 
      + O \left( 
      \left( 1 + |x| X^c \right) 
       \,  \max_{y \le X} \, \max_{(l, d)=1} \, |\Delta(y, d, l)|
      \right)
      \notag \\
      & \notag \\
      & =
      \frac{I(x)}{\varphi(d)} + 
       O \left( 
       (\log X)^B
       \,  \max_{y \le X} \, \max_{(l, d)=1} \, |\Delta(y, d, l)|
      \right) .
      \label{379}
\end{align}
It remains to substitute the above expression for $S(x, d)$ in the formula for 
$L_1(x)$ in \eqref{340}. Since $\delta < \frac{1}{2}$, we are in position to apply 
Lemma~\ref{lemma4}
and find that the contribution from the remainder term 
in \eqref{379} is $\ll X (\log X)^{-A}$.
We also take into account \eqref{330} to conclude
that in the case \eqref{350} the asymptotic formula \eqref{300} is true.

\bigskip

Consider now the case
\begin{equation} \label{380}
  X^{-c} (\log X)^B < |x| \le \tau .
\end{equation}
We proceed with the sum $S(x, d)$, specified by \eqref{340}, in a different way.
Using the properties of Dirichlet's characters we get
\begin{align}
  S(x, d) 
   & = \sum_{\mu X < n \le X} \Lambda(n) \, e \left(x n^c \right) \;
   \frac{1}{\varphi(d)} \sum_{\chi \!\!\!\! \pmod{d} } \chi(n) \, \overline{\chi}(-2)
    \notag \\
  & \notag \\
  & =
   \frac{1}{\varphi(d)}  \sum_{\chi \!\!\!\! \pmod{d} } \overline{\chi}(-2)
     \sum_{\mu X < n \le X} \Lambda(n) \, \chi(n) \, e \left(x n^c \right) .
     \notag
\end{align}
We separate the contribition from the principal character
to get
\begin{equation} \label{390}
  S(x, d) = \frac{1}{\varphi(d)}  Y(X) 
  + O \left( \frac{(\log X)^2}{\varphi(d)} \right) + 
  \mathfrak N ,
\end{equation}
where 
\begin{equation} \label{391}
  Y(X) = \sum_{\mu X < n \le X } \Lambda(n) \, e \left(x n^c \right) 
\end{equation}
and where $\mathfrak N $ comes from the non-principal characters.
We express $\mathfrak N $ as a sum over primitive characters and find that
\[
   \mathfrak N = \frac{1}{\varphi(d)} \sum_{\substack{r \mid d \\ r>1 }} \;
   \sum_{\chi \!\!\!\! \pmod{r}^* } \overline{\chi}(-2) \, 
    \sum_{\substack{\mu X < n \le X \\ (n, d)=1 }} \Lambda(n) \, \chi(n) \, e \left(x n^c \right) .
\]
If we omit the condition $(n, d)=1$ imposed in the sum over $n$, then the resulting error will be
$O \left( (\log X)^2 \right)$. (We leave the easy verification to the reader).
We use \eqref{380} and substitute the expression for $S(x, d)$ in the sum $L_1(x)$
from \eqref{340}. Taking into account \eqref{330} we find that
\begin{equation} \label{400}
  L(x) = \sum_{d \le D} \frac{\lambda(d)}{\varphi(d)} \; 
  Y(X)
  + O \left( X^{\frac{1}{2} + \varepsilon }\right) 
  + \mathfrak M ,
\end{equation}
where
\[
  \mathfrak M = \sum_{d \le D} \frac{\lambda(d)}{\varphi(d)} \;
  \sum_{\substack{r \mid d \\ r>1 }} \;
   \sum_{\chi \!\!\!\! \pmod{r}^* } \overline{\chi}(-2) \, 
    Y(X, \chi)  ,
\]
and
\begin{equation} \label{410}
    Y(X, \chi) = \sum_{\mu X < n \le X } \Lambda(n) \, \chi(n) \, e \left(x n^c \right) .
\end{equation}

\bigskip

We change the order of summation to get
\[
  \mathfrak M = \sum_{1 < r \le D } \;
  \Big( \, \sum_{\substack{d \le D \\ r \mid d}} \frac{\lambda(d)}{\varphi(d)} \, \Big) \;
    \sum_{\chi \!\!\!\! \pmod{r}^* } \overline{\chi}(-2) \, 
    Y(X, \chi) .
\]
Using \eqref{288}, we easily find that
\[
  \sum_{\substack{d \le D \\ r \mid d}} \frac{\lambda(d)}{\varphi(d)} \ll \frac{\log X}{\varphi(r)} ,
\]
hence
\begin{equation} \label{420}
  \mathfrak M \ll (\log X) \sum_{1 < r \le D } \frac{1}{\varphi(r)}
    \sum_{\chi \!\!\!\! \pmod{r}^* }  \, 
    \left| Y(X, \chi)  \right| .
\end{equation}

\bigskip

Consider the sum $Y(X, \chi)$. We apply Abel's formula and use \eqref{55.1} to find that
\begin{align}
  Y(X, \chi) 
    & =
    e \left(x X^c \right) \sum_{\mu X < n \le X } \Lambda(n) \, \chi(n) 
   -
   \int_{\mu X}^X 
   \Big(
   \sum_{\mu X < n \le t} \Lambda(n) \, \chi(n) 
   \Big)
   \; \frac{d}{d t} e \left( x t^c \right) 
   \, d t
   \notag \\
   & \notag \\
   & =
       e \left(x X^c \right) \left( \psi(X, \chi) - \psi(\mu X, \chi) \right)
       -
       \int_{\mu X}^X
       \left( \psi(t, \chi) - \psi(\mu X, \chi) \right)
       \; \frac{d}{d t} e \left( x t^c \right) 
   \, d t .
      \notag
\end{align}
We choose 
\begin{equation} \label{465}
  T = |x| X^c D (\log X)^B .
\end{equation}
From \eqref{102}, \eqref{103} and \eqref{380} we easily see that $2 \le T \le \mu X$.
Now we apply  formula \eqref{450} from Lemma~\ref{lemma4.5} and
find that
\begin{align}
  Y(X, \chi)
   & =
    e \left( x X^c\right) \,
     \left( 
      - \sum_{|\gamma| < T} \frac{X^{\rho} - (\mu X)^{\rho}}{\rho} + O \left( \frac{X (\log X)^2}{T} \right)
    \right) 
   \notag \\
   & \notag \\
  & \qquad \qquad \qquad +
  \int_{\mu X}^X 
    \left(
      \sum_{|\gamma| < T} \frac{t^{\rho} - (\mu X)^{\rho}}{\rho} 
      + O \left( \frac{X (\log X)^2}{T} \right)
    \right) 
    \; \frac{d}{d t} e \left( x t^c \right) \, d t .
  \notag
\end{align}
We estimate the contributions form the error terms, 
then we change the order of summation and integration 
and finally we integrate by parts to get
\[
  Y(X, \chi)
   = - \sum_{|\gamma| < T} I_{\rho} (x) + O \left( \frac{X}{T} (\log X)^2 \left( 1 + |x| X^c \right) \right) ,
\]
where
\begin{equation} \label{500}
  I_{\rho} (x) = \int_{\mu X}^X t^{\rho - 1} e \left( x t^c \right) \, d t .
\end{equation}

\bigskip

Now, we substitute the last expression for $ Y(X, \chi)$ in \eqref{420} 
and use \eqref{380} 
to obtain
\begin{equation} \label{520}
 \mathfrak M \ll (\log X) \sum_{1 < r \le D} \frac{1}{\varphi(r)} \; 
 \sum_{\chi \!\!\!\! \pmod{r}^* } \; \sum_{|\gamma| \le T} \left| I_{\rho} (x) \right|
  \; + \; T^{-1} \, |x| \, X^{1 + c} \, D \, (\log X)^3 .
\end{equation}

\bigskip

We study the sum $Y(X)$, defined by \eqref{391}, in the same manner but now we apply 
formula \eqref{449} from Lemma~\ref{lemma4.5}.
We take the parameter $T$, defined by \eqref{465}, and
after some calculations we find that
\[
   Y(X) = I(x) - \sum_{|\gamma| < T} I_{\rho} (x) + O \left( T^{-1} |x| X^{1 + c} (\log X)^3 \right) ,
\]
where $I(x)$ and $I_{\rho}(x)$ are defined respectively by \eqref{290} and \eqref{500}
and $\rho = \beta + i \gamma $ runs over the non-trivial zeros of $\zeta(s)$ such that $|\gamma| < T$.

\bigskip

Substituting the last expression for $Y(X)$ in \eqref{400}, together with \eqref{465} and \eqref{520}
gives
\begin{equation} \label{540}
 L(x) - \sum_{d \le D} \frac{\lambda(d)}{\varphi(d)} \; I(x) 
 \ll
      X(\log X)^{3 - B}   + (\log X)^3 \mathfrak K  ,
\end{equation}
where
\begin{equation} \label{560}
  \mathfrak K = \sum_{d \le D} \frac{1}{d} \, \sum_{\chi \!\!\!\! \pmod{d}^*} \;
  \sum_{|\gamma| \le T } \left| I_{\rho} (x)\right| .
\end{equation}

\bigskip

Consider the integral $I_{\rho}(x)$ defined by \eqref{500} and let $\rho = \beta + i \gamma$.
With a change of variables we can write the integral in the form
\[
  I_{\rho}(x) = \frac{1}{c} \int_{(\mu X)^c}^{X^c} 
    u^{\frac{1}{c} \rho - 1} e( x u ) \, d u 
    = \int_{(\mu X)^c}^{X^c} u^{\frac{1}{c} \beta - 1} e ( h(u)) \, d u ,
\]
where
\[
  h(u) = \frac{\gamma}{2 \pi c} \log u + x u .
\]
We have
\[
  h'(u) = \frac{\gamma}{2 \pi c \, u}  + x , \qquad
  h''(u) = \frac{- \gamma}{2 \pi c \, u^2} .
\]

\bigskip

Suppose that $|\gamma| \ge 4 \pi c |x| X^c$. Then we have $|h'(u)| \asymp \frac{|\gamma|}{X^c}$ and we estimate the integral using Lemma~\ref{lemma3.5}.
We find that
\begin{equation} \label{580}
   I_{\rho}(x) \ll \frac{X^{\beta}}{|\gamma|}  \qquad \text{if} \qquad 4 \pi c |x| X^c \le |\gamma| \le T .
\end{equation}

\bigskip

If $\pi c \mu^c |x| X^c < |\gamma| < 4 \pi c |x| X^c $,
then we have $|h''(u)| \asymp \frac{|\gamma|}{X^{2c}} \asymp \frac{|x|}{X^c}$ and using
Lemma~\ref{lemma3.5} we find that
\[
  I_{\rho}(x) \ll \frac{X^{\beta - c}}{\sqrt{|x| X^{-c}}} = \frac{X^{\beta}}{\sqrt{|x| X^c}} .
\]

\bigskip

Finally, if $|\gamma| \le \pi c \mu^c |x| X^c$, then we have $|h'(u)| \asymp |x|$ and 
applying again Lemma~\ref{lemma3.5} we get
\[
   I_{\rho} (x) \ll \frac{X^{\beta}}{|x| X^c} .
\]
However, from \eqref{380} it follows that $|x| X^c \ge \sqrt{|x| X^c}$, hence from the last two estimates
we obtain
\begin{equation} \label{600}
  I_{\rho}(x) \ll \frac{X^\beta}{\sqrt{|x| X^c}} \qquad \text{if} \qquad |\gamma| < 4 \pi c |x| X^c .
\end{equation}

\bigskip

From \eqref{560} -- \eqref{600} we find that
\begin{equation} \label{620}
  \mathfrak K \ll \mathfrak K' + \mathfrak K'' ,
\end{equation}
where $\mathfrak K'$ is the contribution comming from the terms for which the condition \eqref{580} for $\gamma$
is satisfied and respectively, 
$\mathfrak K''$ comes from the terms for which $\gamma$ satisfies the condition for $\gamma$ given in \eqref{600}.

\bigskip

Consier $\mathfrak K'$. We have
\[
  \mathfrak K' = \sum_{d \le D} \frac{1}{d} \, \sum_{\chi \!\!\!\! \pmod{d}^*} \;
  \sum_{4 \pi c |x| X^c \le |\gamma| \le T} \frac{X^{\beta}}{|\gamma|} .  
\]
We divide the sum over $d$ into $O (\log X)$ sums in which $d$ runs over an interval of the form
$(Q, Q']$, where $Q' \le \min(2 Q, D) $. We proceed with the sum over $\gamma$ in the same way.
Hence, we obtain
\begin{equation} \label{640}
  \mathfrak K' \ll (\log X)^2 \, \max_{Q \in [1, D]} \; \max_{L \in [4 \pi c |x| X^c, T]} \, 
  \left( \, (QL)^{-1} \, \mathfrak L (L, Q, X) \, \right) ,
\end{equation}
where $\mathfrak L (L, Q, X)$ is the sum defined by \eqref{660}.

\bigskip

Having in mind \eqref{102}, \eqref{103}, \eqref{380}, \eqref{465} and 
since for the above sums we have $Q \le D$ and $L \le T$, we easily verify the condition 
$Q^2 L \le X^{\frac{12}{25}}$ imposed in Lemma~\ref{lemma5}.
Then using \eqref{680}, \eqref{380} and \eqref{640}, we find that
\begin{equation} \label{700}
  \mathfrak K' \ll X (\log X)^{22 - \frac{B}{2}} .
\end{equation}

\bigskip

Consider now the quantity $\mathfrak K''$ for which we have
\[
 \mathfrak K'' \ll \frac{1}{\sqrt{|x| X^c}} \sum_{d \le D}  \frac{1}{d} \;
   \sum_{\chi \!\!\!\! \pmod{d}^*} \; \sum_{|\gamma| \le 4 \pi c |x| X^c} X^{\beta} .
\]
Using \eqref{660}, the estimate \eqref{680} from Lemma~\ref{lemma5} and \eqref{380}  we find that
\begin{equation} \label{720}
   \mathfrak K'' \ll  \frac{\log X}{\sqrt{|x| X^c}} \max_{Q \in [1, D]} 
   \left( \, Q^{-1} \, \mathfrak L \left(4 \pi c |x| X^c, Q, X \right) \, \right)
   \ll X (\log X)^{22 - \frac{B}{2}} .
\end{equation}

\bigskip

From \eqref{620}, \eqref{700} and \eqref{720} we obtain
\[
 \mathfrak K \ll  X (\log X)^{22 - \frac{B}{2}} .
\]
We substitute this estimate for $\mathfrak K$ in \eqref{540} and, as the constant $B$ can be taken arbitrariy large, 
we see that the asymptotic formula \eqref{300} is correct also in the case \eqref{380}.
This proves the lemma.

\hfill
$\square$

\bigskip

The next lemma is an analog of \cite[Lemma~7]{Tol_1}.
\begin{lemma} \label{lemma30}
Let $\tau$ and $\Xi$ be defined by \eqref{103} and \eqref{103.5}.
Then for the sum $L(x)$ and for the integral $I(x)$, defined respectively 
by \eqref{289} and \eqref{290}, we have
\begin{align}
   \int_{|x| < \tau} \left| L(x) \right|^2 \, d x 
    & \ll 
    X^{2-c} (\log X)^6 , \label{740} \\
    & \notag \\
  \int_{|x| < \tau} \left| I(x) \right|^2 \, d x 
   & \ll 
   X^{2-c} (\log X)^4 , \label{742} \\
   & \notag \\
   \int_{|x| <   \Xi } \left| L(x) \right|^2 \, d x 
    & \ll 
    X \, \Xi \, (\log X)^6 . \label{744} 
\end{align}
\end{lemma}

{\bf Proof:}
We shall only prove \eqref{740}, the other inequalities can be proved likewise.
Denote by $J$ the integral on the left side of \eqref{740}. We have
\begin{align}
 J
  & = 
  \int_{|x| < \tau} L(x) L(-x) \, d x
  \notag \\
  & \notag \\
  & =
    \int_{|x| < \tau}
     \sum_{d_1 \le D } \lambda(d_1) \sum_{\substack{\mu X < p_1 \le X \\ d_1 \mid p_1 + 2 }}
       (\log p_1) e \left( x p_1^c \right) \; 
        \sum_{d_2 \le D } \lambda(d_2) \sum_{\substack{\mu X < p_2 \le X \\ d_2 \mid p_2 + 2 }}
       (\log p_2) e \left( - x p_2^c \right) \, d x
     \notag \\
  & \notag \\
  & =
   \sum_{d_1, d_2 \le D} \lambda(d_1) \lambda(d_2) 
     \sum_{\substack{\mu X < p_1, p_2 \le X \\ d_i \mid p_i + 2 \; , \; i=1,2}}
     (\log p_1 ) (\log p_2) \int_{|x| < \tau} 
      e \left( x \left( p_1^c - p_2^c\right) \right) \, d x
       \notag \\
  & \notag \\
  & \ll
   (\log X)^2
    \sum_{d_1, d_2 \le D} \;
     \sum_{\substack{\mu X < n_1, n_2 \le X \\ d_i \mid n_i + 2 \; , \; i=1,2}}
      \min \left( \tau, \frac{1}{|n_1^c - n_2^c|} \right) .
   \notag
\end{align}
We change the order of summation and use the obvious inequality $uv \le u^2 + v^2$ to get
\begin{align}
 J 
 & \ll
  (\log X)^2
   \sum_{\mu X < n_1, n_2 \le X }
    \tau(n_1 + 2 ) \, \tau(n_2 + 2) \, 
      \min \left( \tau, \frac{1}{|n_1^c - n_2^c|} \right) 
        \notag \\
  & \notag \\
  & \ll
  (\log X)^2
   \sum_{\mu X < n_1, n_2 \le X }
    \tau^2(n_1 + 2 )  \, 
      \min \left( \tau, \frac{1}{|n_1^c - n_2^c|} \right) .
 \notag 
\end{align}
Now we proceed as in the proof of \cite[Lemma~7]{Tol_1} and we also
use the well-known inequality $\sum_{n \le y} \tau^2 (n) \ll y (\log y)^3$.
In this way we prove \eqref{740} ---
we leave the details to the reader.

\hfill
$\square$

\bigskip

We shall find asymptotic formulae for the integrals
$\Gamma_1^{(1)}$ and $\Gamma_4^{(1)}$ defined respectively by \eqref{260} and \eqref{260.1}.
We consier only $\Gamma_1^{(1)}$ because the study of $\Gamma_4^{(1)}$ is similar.

\bigskip

From \eqref{184}, \eqref{290} and  Lemma~\ref{lemma10} we know that if $|x| < \tau $, then we have
\[
  L^{\pm}(x) = \mathfrak N^{\pm} I(x) + O \left( X (\log X)^{-A} \right) ,
\]
where $A>0$ is arbitrarily large.
Let
\begin{equation} \label{745}
  \mathfrak I_0 = \int_{|x|< \tau} \Theta(x) \, e \left( - N x \right) \, I(x)^3 \, d x .
\end{equation}
We use the identity
\begin{align}
  &
  L^{-}(x) \left( L^{+}(x) \right)^2 
   = \mathfrak N^{-} \left( \mathfrak N^{+} \right)^2 I(x)^3
   + 
  \left(  
   L^{-}(x)  - \mathfrak N^{-} I(x) \right) \left( \mathfrak N^{+} \right)^2 I(x)^2
  \notag \\
  & \notag \\
    & \qquad \qquad  +
     L^{-}(x)  \left( L^{+}(x)  - \mathfrak N^{+} I(x) \right) \mathfrak N^{+}  I(x)
    \; + \; 
    L^{-}(x)   L^{+}(x) \left( L^{+}(x)  - \mathfrak N^{+} I(x) \right)
  \notag
\end{align}
and the obvious estimate 
\begin{equation} \label{848}
 \mathfrak N^{\pm} \ll \log X
\end{equation}
to find that
\[
   \left| L^{-}(x) \left( L^{+}(x) \right)^2 
   - \mathfrak N^{-} \left( \mathfrak N^{+} \right)^2 I(x)^3 \right|
   \ll 
   X (\log X)^{3-A} \left( \left| I(x) \right|^2 +\left| L^- (x) \right|^2 + \left| L^+ (x) \right|^2 \right) .
\]
From the above inequality, the estimate \eqref{158} with $a = \frac{7}{8} \Delta$, from 
\eqref{260}, \eqref{745} and
the estimates \eqref{740}, \eqref{742} from Lemma~\ref{lemma30}
it follows that
\begin{equation} \label{849}
  \Gamma_1^{(1)} -  \mathfrak N^{-} \left( \mathfrak N^{+} \right)^2 \, \mathfrak I_0
  \ll \Delta X^{3-c} (\log X)^{9 - A} .
\end{equation}

\bigskip

Consider now the integral
\begin{equation} \label{850.5}
   \mathfrak I = \int_{- \infty}^{\infty} \Theta(x) \, e \left( - N x\right) \, I(x)^3 \, d x .
\end{equation}
Following the proof of 
\cite[Lemma~6]{Tol_1} we see that for a suitable $\mu \in (0, 1)$ depending on $c$
we have
\begin{equation} \label{850.6}
   \mathfrak I  \gg 
    \Delta X^{3 - c} .
\end{equation}
We leave the easy verification to the reader.

\bigskip

Further, we have
$\left| \frac{d}{d t} (x t^c)\right| \gg |x| X^{c-1}$
for $t \in [\mu X, X]$. Therefore,
applying again Lemma~\ref{lemma3.5}
we find that
$I(x) \ll |x|^{-1} X^{1 - c}$.
From this estimate, \eqref{103}, \eqref{158} with $a = \frac{7}{8} \Delta$, \eqref{745} and \eqref{850.5}
we find that
\begin{equation} \label{850.7}
  |\mathfrak I - \mathfrak I_0| \le \int_{|x| > \tau} |\Theta(x)| \, |I(x)|^3 \, d x
  \ll \Delta  \, X^{3 - 3c } \int_{x > \tau} \frac{d x}{x^3}
  \ll  \Delta \, X^{3 - 3c }  \, \tau^{-2} 
  \ll \Delta X^{3 - c - 2 \xi} .
\end{equation}

\bigskip

We apply \eqref{849} (with $A = 12$) as well as \eqref{848} and \eqref{850.7}
to find
\begin{equation} \label{850.8}
  \Gamma_1^{(1)} = \mathfrak N^{-} \left( \mathfrak N^{+} \right)^2 \mathfrak I
  + O \left( \Delta X^{3-c} (\log X)^{-4} \right) .
\end{equation}
We proceed with $\Gamma_4^{(1)}$ is the same way and prove that
\begin{equation} \label{850.9}
  \Gamma_4^{(1)} = \left( \mathfrak N^{+} \right)^3 \mathfrak I
  + O \left( \Delta X^{3-c} (\log X)^{-4} \right) .
\end{equation}

\section{The estimation of $\Gamma_1^{(2)}$ and $\Gamma_4^{(2)}$ and the end of the proof}

From \eqref{288.5} we see that in order to find a non-trivial lower bound for $\Gamma$ we have to prove that
the integrals $\Gamma_1^{(2)}$ and $\Gamma_4^{(2)}$ are small enough. To establish this 
we need estimates for the sums $L^{\pm}(x)$ provided that $\tau \le |x| \le \Xi$. 

\bigskip

We apply the next lemma, which is a special case of Vaughan's identity.
\begin{lemma} \label{lemma35}
Let $f(n)$ be a complex valued function defined for integers $n \in \left(\mu X, X \right]$.
Then we have
\[
 \sum_{\mu X < n \le X} \Lambda(n) f(n) = S_1 - S_2  - S_3 ,
\]
where
\begin{align}
   S_1 
     & =  \sum_{k \le X^{\frac{1}{3}}} \mu(k) \sum_{\frac{\mu X}{k} < l \le \frac{X}{k}} (\log l) f(kl) ,
  \label{920} \\
  & \notag \\
   S_2
     & = 
  \sum_{ k \le X^{\frac{2}{3}}} c(k) \sum_{\frac{\mu X}{k} < l \le \frac{X}{k}}   f(kl) ,
    \label{940} \\
  & \notag \\
   S_3
     & = 
     \sum_{ X^{\frac{1}{3}} < k \le X^{\frac{2}{3}}} a(k) 
       \sum_{\frac{\mu X}{k} < l \le \frac{X}{k}}  \Lambda(l) f(kl) 
     \label{960}
\end{align}
and where $a(k)$, $c(k)$ are real numbers satisfying
\begin{equation} \label{980}
   |a(k)| \le \tau(k) ,  \qquad  |c(k)| \le \log k .
\end{equation}
\end{lemma}

{\bf Proof.}
Can be found in \cite{Vaughan}.

\hfill
$\square$

\bigskip

Follows Van der Corput's inequality.
\begin{lemma} \label{lemma45}
Let $\alpha$, $\beta$ be 
real numbers with $\beta - \alpha \ge 1$
and let $H$ be a positive integer. 
Suppose that for any integer $l \in (\alpha, \beta ]$
there is a complex number $\Upsilon (l)$.
Then we have
\begin{equation} \label{860}
 \left| \sum_{\alpha < l \le \beta} \Upsilon (l) \right|^2
 \le
 \frac{\beta - \alpha + H}{H} \sum_{|h| < H} 
 \left( 1 - \frac{|h|}{H} \right) \sum_{\alpha < l, l+h \le \beta} \Upsilon (l+h) \,  \overline{\Upsilon (l)} .
\end{equation}
\end{lemma}

{\bf Proof.}
Can be found in \cite[Lemma~8.17]{Iv_Kow}.

\hfill
$\square$

\bigskip

The next lemma presents Van der Corput's estimate for 
exponential sums.
\begin{lemma} \label{lemma55}
Suppose that $\alpha, \beta$ are real numbers with $\beta - \alpha \ge 1$ and let 
$f(y)$ be two times continuously differentiable function in the interval $[\alpha, \beta]$.
Assume also that for some $\lambda > 0$ we have 
$ \left| f''(y)\right| \asymp \lambda $ uniformly for $y \in [\alpha, \beta]$.
Then we have
\[ 
  \left| \sum_{\alpha < n \le \beta} e \left( f(n) \right) \right| 
  \ll (\beta - \alpha) \lambda^{\frac{1}{2}} + \lambda^{- \frac{1}{2}} .
\]
\end{lemma}

{\bf Proof.}
See \cite[Chapter~1, Theorem~5]{Karatsuba}.

\hfill
$\square$

\bigskip

From this point onwards we assume that
\begin{equation} \label{875}
  \xi = \frac{459 c - 435}{125} , \qquad \delta = \frac{180 - 168 c}{125} .
\end{equation}
(It is easy to verify that \eqref{875} implies the first inequality from \eqref{102}).
We prove the following
\begin{lemma} \label{lemma100}
Suppose that $D$, $\tau$ and $ \Xi $ are defined by \eqref{103} and \eqref{103.5} 
and $\xi$, $\delta$ are specified by \eqref{875}.
Suppose also that
the real numbers $\lambda(d)$ satisfy \eqref{288} and $L(x)$ is defined by \eqref{289}.
Then there exists $\varkappa(c) > 0$ such that
\begin{equation} \label{880}
  \sup_{\tau \le |x| \le \, \Xi } \left| L(x) \right| \ll X^{2 - c - \varkappa(c)} .
\end{equation}
\end{lemma}

{\bf Proof.}
Instead of $L(x)$ we consider the sum $L_1(x)$ given by \eqref{340}
and we take into account \eqref{330}.
We write $L_1(x)$ in the form
\begin{equation} \label{895}
  L_1(x) = \sum_{\mu X < n \le X} \Lambda(n) \, f(n) , \qquad
  f(n) = \sum_{\substack{d \le D \\ d \mid n+2}} \lambda(d) \, e \left( x n^c \right)  
\end{equation}
and we apply Lemma~\ref{lemma35} to find that
\begin{equation} \label{900}
   L_1(x) = S_1 - S_2 - S_3 , 
\end{equation}
where $S_1$, $S_2$ and $S_3$ are the sums defined respectively by \eqref{920} -- \eqref{960}
with the function $f(n)$ given by \eqref{895}.

\bigskip

From \eqref{940} we find that
\begin{equation} \label{1000}
  S_2 = S_2'+ S_2'' ,
\end{equation}
where
\begin{align}
  S_2' 
  & =
   \sum_{k \le X^{\frac{1}{3}}} c(k) \sum_{\frac{\mu X}{k} < l \le \frac{X}{k}}
    f(kl) ,
   \label{1020} \\
   & \notag \\
  S_2''
   & = 
   \sum_{X^{\frac{1}{3}} < k \le X^{\frac{2}{3}}} c(k) \sum_{\frac{\mu X}{k} < l \le \frac{X}{k}}
    f(kl) .
   \label{1040} 
\end{align}
Therefore we have
\begin{equation} \label{1060}
  L_1(x) \ll |S_1| + |S_2'| + |S_2''| + |S_3| .
\end{equation}

\bigskip

Consider the sum $S_2'$. We use \eqref{895} and \eqref{1020} and change the order of summation to
write it in the form
\[
 S_2' = \sum_{d \le D} \lambda(d) \sum_{k \le X^{\frac{1}{3}}} c(k)
 \sum_{\substack{\frac{\mu X}{k} < l \le \frac{X}{k} \\ d \mid kl+2 }}
 e \left( x (kl)^c \right) .
\]
Since $\lambda(d)= 0 $ for $2 \mid d$, then from $d \mid kl + 2$ it follows that
$(k, d)=1$. Hence there exists an integer $l_0$ such that $d \mid kl + 2$
is equivalent to $l \equiv l_0 \pmod{d}$, which means that $l = l_0 + m d$ for some integer $m$.
Therefore we get
\begin{equation} \label{1080}
  S_2' = 
  \sum_{d \le D} \lambda(d) \sum_{\substack{k \le X^{\frac{1}{3}} \\ (k, d)=1 } } c(k)
 \sum_{\frac{\mu X}{kd} - \frac{l_0}{d} < m \le \frac{X}{kd} - \frac{l_0}{d} }
 e \left( h(m) \right) ,
 \qquad  h(m) = x k^c (l_0 + md)^c .
\end{equation}
We have $h''(m) = c(c-1)x k^c d^2 (l_0 + md)^{c-2}$, hence
$|h''(m)| \asymp |x| k^2 d^2 X^{c-2}$ and using 
Lemma~\ref{lemma55} we find that the sum over $m$ in \eqref{1080}
is 
\[
 \ll \frac{X}{kd} \left( |x| k^2 d^2 X^{c-2} \right)^{\frac{1}{2}}
 + \left( |x| k^2 d^2 X^{c-2} \right)^{- \frac{1}{2}}
 \ll
  |x|^{\frac{1}{2}} X^{\frac{c}{2}} + |x|^{- \frac{1}{2}} (kd)^{-1} X^{1 - \frac{c}{2}} .
\]
Then using \eqref{288}, \eqref{980} and \eqref{1080}
we find that
\begin{equation} \label{1100}
  S_2' \ll
  X^{\varepsilon} \left( D |x|^{\frac{1}{2}} X^{\frac{1}{3} + \frac{c}{2}} 
  + |x|^{- \frac{1}{2}} X^{1 - \frac{c}{2}} \right) .
\end{equation}

\bigskip

To etimate $S_1$ we apply Abel's summation formula to
get rid of the factor $(\log l)$ in the sum from \eqref{920} and then proceed as above.
In this way we find that
\begin{equation} \label{1120}
  S_1 \ll
  X^{\varepsilon} \left( D |x|^{\frac{1}{2}} X^{\frac{1}{3} + \frac{c}{2}} 
  + |x|^{- \frac{1}{2}} X^{1 - \frac{c}{2}} \right) .
\end{equation}
We leave the simple calculations to the reader.

\bigskip

Consider now the sum $S_3$.
We divide it into $O \left( \log X \right)$ sums of type
\begin{equation} \label{1200}
  W(K) =
   \sum_{K < k \le K_1} a(k) \;
   \sum_{\frac{\mu X}{k} < l \le \frac{X}{k}} \Lambda(l)
   \sum_{\substack{ d \le D \\ d \mid kl+2}} \lambda(d) e \left( x (kl)^c \right) ,
\end{equation}
where
\begin{equation} \label{1220}
   K_1 \le 2 K , \qquad X^{\frac{1}{3}} \le K < K_1 \le X^{\frac{2}{3}} .
\end{equation}

\bigskip

Consider the case
\begin{equation} \label{1230}
   X^{\frac{1}{2} } \le K .
\end{equation}
It is clear that
\begin{equation} \label{1240}
    W(K) 
    \ll 
    X^{\varepsilon}
   \sum_{K < k \le K_1} 
    \left| 
   \sum_{\frac{\mu X}{k} < l \le \frac{X}{k}}  \Upsilon (l)
   \right| ,
\end{equation}
where
\begin{equation} \label{1250}
   \Upsilon  (l) = \Lambda(l) \; \sum_{\substack{ d \le D \\ d \mid kl+2}} \lambda(d) \, \, e \left( x (kl)^c \right) .
\end{equation}
Applying Cauchy's inequality, we find that
\begin{equation} \label{1260}
    |W(K)|^2 
    \ll 
    X^{\varepsilon} K 
   \sum_{K < k \le K_1} 
    \left| 
   \sum_{\frac{\mu X}{k} < l \le \frac{X}{k}}  \Upsilon (l)
   \right|^2 .
\end{equation}

\bigskip

Assume that $H$ is an integer satisfying
\begin{equation} \label{1270}
    1 \le H \ll \frac{X}{K}
\end{equation}
and let $\Upsilon  (l) $ be defined by \eqref{1250}. 
For the inner sum in \eqref{1260} we apply 
the inequaity
\eqref{860} from Lemma~\ref{lemma45} and we find that
\begin{align}
  |W(K)|^2
   & 
   \ll
   \frac{X^{1 + \varepsilon}}{H}
   \sum_{K < k \le K_1} \; \sum_{|h| < H}
     \left( 1 - \frac{|h|}{H} \right) \;
       \sum_{\frac{\mu X}{k} < l, l+h \le \frac{X}{k}}
        \notag \\
        & \notag \\
    & \qquad \qquad  
           \Lambda(l+h) \sum_{\substack{d_2 \le D \\ d_2 \mid k(l+h)+2}}
        \lambda(d_2) \, e \left( x(k(l+h))^c \right) 
        \;\;
          \Lambda(l) \sum_{\substack{d_1 \le D \\ d_1 \mid kl+2}}
        \lambda(d_1) \, e \left( - x(kl)^c \right) .
     \notag
\end{align}
We change the order of summation and find
\begin{equation} \label{1280}
  |W(K)|^2
      \ll
   \frac{X^{1 + \varepsilon}}{H}
   \sum_{d_1, d_2 \le D} \lambda(d_1) \lambda(d_2)
   \sum_{|h| < H}
    \left( 1 - \frac{|h|}{H} \right)
     \sum_{\frac{\mu X}{K_1} < l, l+h \le \frac{X}{K}}
      \Lambda(l) \Lambda(l+h) \; \mathfrak F ,
\end{equation}
where
\[
  \mathfrak F = \sum_{\substack{\widetilde{K} < k \le \widetilde{K_1} \\ d_1 \mid kl+2 \\ d_2 \mid k(l+h) + 2}}
  e \left( x k^c \left( (l+h)^c - l^c \right) \right)
\]
and
\[
   \widetilde{K} = \max \left( K , \frac{\mu X}{l} , \frac{\mu X}{l+h} \right) , \qquad
   \widetilde{K_1} = \min \left( K_1 , \frac{X}{l} , \frac{X}{l+h} \right) .
\]

\bigskip

Since $\lambda(d)= 0 $ for $2 \mid d$, we may assume that $2 \nmid d_1 d_2$. 
Hence $(d_1, l) =(d_2, l+h)=1 $ because otherwise the sum $\mathfrak F$ would be empty.
Therefore, there exists an integer $k_0$ depending on $l, h, d_1, d_2$ such that 
the pair of conditions $ d_1 \mid kl+2 $ and $ d_2 \mid k(l+h) + 2$ is equivalent to 
the congruence $k \equiv k_0 \pmod{[d_1, d_2]}$.
Hence we may write the sum $\mathfrak F$ as
\[
  \mathfrak F = 
  \sum_{\frac{\widetilde{K} - k_0}{[d_1, d_2]} < m \le \frac{\widetilde{K_1}- k_0}{[d_1, d_2]}} 
  e(F(m)) ,
\]
where
\[
  F(m) = x \left( (l+h)^c  - l^c \right) \left( k_0 + m [d_1, d_2] \right)^c .
\]

\bigskip

Obviously, we have
\begin{equation} \label{1290}
  \mathfrak F  \ll \frac{K}{[d_1, d_2]} \qquad \text{if} \qquad h = 0 .
\end{equation}
Consider now the case $h \not= 0$.
We have
\[
  F''(m) = c (c-1) x \left( (l+h)^c  - l^c \right) \left( k_0 + m [d_1, d_2] \right)^{c-2} [d_1, d_2]^2 ,
\]
hence
\[
  |F''(m)| \asymp  |x| \left| (l+h)^c  - l^c \right| K^{c-2} [d_1, d_2]^2 .
\]
We apply Lemma~\ref{lemma55} and find that
\begin{align}
  \mathfrak F  
  & \ll 
  \frac{K}{[d_1, d_2]} \,  \Big( |x| \left| (l+h)^c  - l^c \right| K^{c-2} [d_1, d_2]^2 \Big)^{\frac{1}{2}}
  + 
  \Big( |x| \left| (l+h)^c  - l^c \right| K^{c-2} [d_1, d_2]^2 \Big)^{-\frac{1}{2}} 
   \notag \\
   & \notag \\
  & \ll
    |x|^{\frac{1}{2}} K^{\frac{c}{2}} \left| (l+h)^c - l^c \right|^{\frac{1}{2}}
      + |x|^{- \frac{1}{2}} K^{1 - \frac{c}{2}} \left| (l+h)^c - l^c \right|^{-\frac{1}{2}}
      [d_1, d_2]^{-1} .
  \notag 
\end{align}
We use the above estimate, \eqref{1280}, \eqref{1290} as well as the estimate 
\[
  \sum_{d_1, d_2 \le D} \frac{1}{[d_1, d_2]} \ll (\log X)^3 
\]
(we leave the easy verification to the reader) to obtain
\begin{align}
  |W(K)|^2
    \ll
    X^{2 + \varepsilon} H^{-1}
    + 
    X^{1 + \varepsilon} \, D^2  \, |x|^{\frac{1}{2}}  \, K^{\frac{c}{2}}  \,  H^{-1} \,
        \Sigma_1 
    +
    X^{1 + \varepsilon} \,  |x|^{-\frac{1}{2}}  \, K^{1 - \frac{c}{2}}  \,  H^{-1} \, 
      \Sigma_1 ,
     \notag
\end{align}
where
\begin{align}
  \Sigma_1 
   & =
    \sum_{0 < |h| \le H} \;
    \sum_{\frac{\mu X}{K_1} < l, l+h \le \frac{X}{K}} 
      \left| (l+h)^c - l^c \right|^{\frac{1}{2}} ,
    \notag \\
    & \notag \\
     \Sigma_2 
     & = 
     \sum_{0 < |h| \le H} \;
    \sum_{\frac{\mu X}{K_1} < l, l+h \le \frac{X}{K}} 
      \left| (l+h)^c - l^c \right|^{- \frac{1}{2}} .
    \notag  
\end{align}
By a straightforward calculation, which we leave to the reader, one obtains
\[
  \Sigma_1 \ll 
  H^{\frac{3}{2}} \, X^{\frac{c+1}{2}} \, K^{- \frac{c+1}{2}} ,
  \qquad
  \Sigma_2
  \ll 
  H^{\frac{1}{2}} \, X^{\frac{3-c}{2}} \, K^{- \frac{3-c}{2}} .
\]
Therefore, we find
\begin{equation} \label{1300}
  |W(K)|^2 
  \ll X^{\varepsilon}
  \left( 
    X^2 H^{-1} + X^{\frac{3 + c}{2}} \, D^2 \, |x|^{\frac{1}{2}} \, K^{- \frac{1}{2}} \, H^{\frac{1}{2}} 
    +
    X^{\frac{5 - c}{2}} \, |x|^{-\frac{1}{2}} \, K^{- \frac{1}{2}} \, H^{-\frac{1}{2}} 
  \right) .
\end{equation}

\bigskip

We choose
\begin{equation} \label{1310}
  H = \left[ \min (H_0, X K^{-1}) \right] , \qquad \text{where} \qquad 
  H_0 = X^{\frac{1-c}{3}} \, K^{\frac{1}{3}} \, D^{-\frac{4}{3}} \, |x|^{- \frac{1}{3}} .
\end{equation}
(It is easy to verify that \eqref{1270} holds).
We note that 
\begin{equation} \label{1320}
 H^{-1} \asymp H_0^{-1} + K X^{-1} .
\end{equation}
Using \eqref{1220}, \eqref{1230} and \eqref{1300} -- \eqref{1320} we obtain
\begin{equation} \label{1330}
  W(K) \ll X^{\varepsilon}
  \left( 
   X^{\frac{3}{4} + \frac{c}{6}} D^{\frac{2}{3}} |x|^{\frac{1}{6}}
   + X^{\frac{5}{6}} 
   + X^{1 - \frac{c}{6}} D^{\frac{1}{3}} |x|^{- \frac{1}{6}} 
   + X^{1 - \frac{c}{4}} |x|^{- \frac{1}{4}} 
  \right) .
\end{equation}

\bigskip

Consider now the sum $W(K)$ defined by \eqref{1200} in the case
\begin{equation} \label{1340}
   K < X^{\frac{1}{2}} .
\end{equation}
We write it in the form
\[
 W(K) =
 \sum_{\frac{\mu X}{K_1} < l \le \frac{X}{K}} \Lambda(l) \;
   \sum_{\max \left( K, \frac{\mu X}{l} \right) < k \le \min \left( K_1, \frac{X}{l} \right) } a(k) \;
   \sum_{\substack{d \le D \\ d \mid kl+2}} \lambda(d) \, e \left( x (kl)^c \right) .
\]
Now we have $\frac{X}{K} \gg X^{\frac{1}{2}}$ and we may proceed as above but with 
r\^oles of $k$ and $l$ reversed. Finally, we estabish again the estimate \eqref{1330}.

\bigskip

Since the sum $S_3$ consists of $O \left( \log X\right)$ sums of type $W(K)$,
it can be estimated by the expression from the right side or \eqref{1330}, too.
We study the sum $S_2''$ in the same manner and we obtain
\begin{equation} \label{1350}
  S_2'', S_3 \ll
  X^{\varepsilon}
  \left( 
   X^{\frac{3}{4} + \frac{c}{6}} D^{\frac{2}{3}} |x|^{\frac{1}{6}}
   + X^{\frac{5}{6}} 
   + X^{1 - \frac{c}{6}} D^{\frac{1}{3}} |x|^{- \frac{1}{6}} 
   + X^{1 - \frac{c}{4}} |x|^{- \frac{1}{4}} 
  \right) .
\end{equation}

\bigskip

From \eqref{330},  \eqref{1060},  \eqref{1100},  \eqref{1120} and  \eqref{1350}
we find that
\[
  L(x) \ll 
  X^{\varepsilon}
  \left( 
  X^{\frac{1}{3} + \frac{c}{2}}   D |x|^{\frac{1}{2}}
  +  X^{1 - \frac{c}{2}} |x|^{- \frac{1}{2}}
  +  X^{\frac{3}{4} + \frac{c}{6}} D^{\frac{2}{3}} |x|^{\frac{1}{6}}
   + X^{\frac{5}{6}} 
   + X^{1 - \frac{c}{6}} D^{\frac{1}{3}} |x|^{- \frac{1}{6}} 
   + X^{1 - \frac{c}{4}} |x|^{- \frac{1}{4}} 
   \right) .
\]

\bigskip

Now we use \eqref{103} to find that if $\tau \le |x| \le \Xi $, then we have
$  X^{\xi - c}  \le |x| \ll X^{\varepsilon} $.
Hence
\[
 L(x) \ll 
 X^{\varepsilon}
  \left( 
  X^{\frac{1}{3} + \frac{c}{2} + \delta}  
   +  X^{\frac{3}{4} + \frac{c}{6} + \frac{2 \delta}{3}} 
   + X^{1 - \frac{\xi}{6} + \frac{\delta}{3}} 
   + X^{1 - \frac{\xi}{4}} 
   \right) .
\]
It remains to use \eqref{101} and \eqref{875} and after a simple calculation we obtain \eqref{880}.

\hfill
$\square$

\bigskip

We are now in position to estimate 
the quantities
$\Gamma_1^{(2)}$ and $\Gamma_4^{(2)}$ defined respectively by 
\eqref{270} and \eqref{270.1}.
We apply Lemma~\ref{lemma100} and use \eqref{158} 
with $a = \frac{7}{8} \Delta$ 
to find that
\[
 \Gamma_1^{(2)}, \Gamma_4^{(2)} \ll
 X^{2 - c - \kappa(c)} \, \Delta \, \int_{|x| \le \Xi } \left| L^+ (x)\right|^2 \, d x .
\]
From the above formula,
the definitions of $\Delta$ and $ \Xi $ 
(see \eqref{102.5} and \eqref{103.5})
and formula \eqref{744} of Lemma~\ref{lemma30} we conclude that
\[
 \Gamma_1^{(2)}, \Gamma_4^{(2)} \ll
  \Delta  X^{3 - c } (\log X)^{-4} .
\]

\bigskip

From the last formula and \eqref{288.5}, \eqref{850.6}, \eqref{850.8}, \eqref{850.9} we conclude that
\begin{equation} \label{1360}
   \Gamma \ge \left| 3 \, \mathfrak N^{-} - 2 \, \mathfrak N^{+} \right|  \left( \mathfrak N^{+} \right)^2 \mathfrak J 
   + O \left( \Delta X^{3 - c} (\log X)^{-4} \right) .
\end{equation}
Now we shall
find a lower bound for the difference $3 \, \mathfrak N^{-} - 2 \, \mathfrak N^{+}$.
It is clear that for the quantity $\mathfrak P$ defined by \eqref{184} we have
\begin{equation} \label{1400}
  \mathfrak P \asymp (\log X)^{-1} .
\end{equation}
From \eqref{181} and \eqref{182} we see that
\begin{equation} \label{1410}
  3 \mathfrak  N^{-} - 2 \mathfrak  N^{+} \ge \mathfrak P \, 
  \left( 3 f(s_0) - 2 F(s_0) \right) + O \left( (\log X)^{- \frac{1}{3} } \right) ,
\end{equation}
where $s_0$ is defined by \eqref{184} and $F(s)$, $f(s)$ are the functions specified by \eqref{185}.
We take $  s_0 = 2,95 $ which means that
\[
  \eta = \frac{\delta}{2,95} = \frac{180 - 168 c}{368,75} .
\]
Hence, using \eqref{185} we see that $3 f(s_0) - 2 F(s_0) > 0$.

\bigskip

It remains to take into account
the lower bound for $\mathfrak N^+$ in \eqref{181}
as well as \eqref{850.6}, \eqref{1400}, \eqref{1410}
to obtain
\[
 \Gamma \gg \Delta X^{3 - c} (\log X)^{-3} .
\]
Therefore $\Gamma > 0$ and
the equation \eqref{100} with $\Delta$ specified by \eqref{102.5}
has a solution in primes $p_1, p_2, p_3$ such that
each of the numbers $p_1 + 2, p_2 + 2, p_3 + 2$ has at most 
$\left[ \frac{369}{180 - 168 c} \right]$ prime factors.
This completes the proof of Theorem~\ref{theorem10}.

\hfill
$\square$

\bigskip
\bigskip

\vbox{
\hbox{Faculty of Mathematics and Informatics}
\hbox{Sofia University ``St. Kl. Ohridsky''}
\hbox{5 J. Bourchier, 1164 Sofia, Bulgaria}
\hbox{ }
\hbox{dtolev@fmi.uni-sofia.bg}}

\end{document}